\documentclass{gtart}

\def\ifplaintex{\expandafter\ifx\csname documentclass\endcsname\relax}

\def\gtp{{\mathsurround=0pt\it $\cal G\mskip-2mu$eometry \&\ 
$\cal T\!\!$opology $\cal P\!$ublications}}  

\def\recd{{\small Received:\qua\receiveddate\ifx\reviseddate\relax
\else\qquad Revised:\qua\reviseddate\fi\par}} 

\def\lognumber#1{\def\thelognumber{#1}}
\def\volumenumber#1{\def\thevolumenumber{#1}}
\def\volumeyear#1{\def\thevolumeyear{#1}}
\def\papernumber#1{\def\thepapernumber{#1}}
\def\pagenumbers#1#2{\def\startpage{#1}\def\finishpage{#2}}
\def\published#1{\def\publishdate{#1}}

\def\received#1{\def\receiveddate{#1}}
\def\revised#1{\def\reviseddate{#1}}
\def\accepted#1{\def\accepteddate{#1}}

\def\asciiauthors#1{\def\theasciiauthors{#1}}
\def\asciiaddress#1{\def\theasciiaddress{#1}}

\def\coverauthors#1{\def\thecoverauthors{#1}}
\long\def\asciiabstract#1{\long\def\theasciiabstract{#1}}

\let\\\par\let\thelognumber\relax\let\thevolumenumber\relax
\let\thepapernumber\relax\let\thevolumeyear\relax\let\startpage\relax
\let\finishpage\relax\let\publishdate\relax\let\receiveddate\relax
\let\reviseddate\relax\let\accepteddate\relax\let\theasciititle\relax
\let\theasciiauthors\relax\let\theasciiaddress\relax
\let\theasciiabstract\relax

\let\thecoverauthors\relax\let\theasciiemail\relax

\ifplaintex
\font\logobig=cmssbx10 scaled 3836
\font\logomed=cmssbx10 scaled 2557
\else
\font\logobig=cmssbx10 scaled 4200
\font\logomed=cmssbx10 scaled 2800
\fi

\long\def\makeagttitle{   
\count0=\startpage
\agt\hfill      
\hbox to 45truept{\vbox to 0pt{\vglue -13truept{\logomed A\kern -.37em{\logobig 
T}\kern -.38em G}\vss}\hss}
\break
{\small Volume \thevolumenumber\ (\thevolumeyear)
\startpage--\finishpage\nl
Published: \publishdate}

\vglue .25truein

{\parskip=0pt\leftskip 0pt plus
1fil\def\\{\par\smallskip}{\Large\bf\thetitle}\par\medskip} \vglue
0.05truein

{\parskip=0pt\leftskip 0pt plus 1fil\def\\{\par}{\sc\theauthors}
\par\medskip}%
 
\vglue 0.03truein 

{\small\leftskip 25truept\rightskip 25truept{\bf Abstract}\stdspace\theabstract

{\bf AMS Classification}\stdspace\theprimaryclass
\ifx\thesecondaryclass\relax\else; \thesecondaryclass\fi\par
{\bf Keywords}\stdspace \thekeywords\par}\vglue 7truept

}   

\ifplaintex
\font\phead=cmsl9 scaled 950
\font\pnum=cmbx10 scaled 913
\font\pfoot=cmsl9 scaled 950
\headline{\vbox to 0pt{\vskip -4.5mm\line{\small\phead\ifnum
\count0=\startpage ISSN 1472-2739 (on-line) 1472-2747 (printed)
\hfill {\pnum\folio}\else\ifodd\count0\def\\{ }%
\ifx\theshorttitle\relax\thetitle\else\theshorttitle\fi\hfill{\pnum\folio}
\else\def\\{ and }{\pnum\folio}\hfill\ifx\theshortauthors\relax\theauthors
\else\theshortauthors\fi\fi\fi}\vss}}
\footline{\vbox to 0pt{\vglue 0mm\line{\small\pfoot\ifnum\count0=\startpage
\copyright\ \gtp\hfill\else
\agt, Volume \thevolumenumber\ (\thevolumeyear)\hfill\fi}\vss}}
\else
\font=cmsl9 scaled 1050
\font\lnum=cmbx10 
\font=cmsl9 scaled 1050
\makeatletter
\def\@oddhead{{\small\lhead\ifnum\count0=\startpage ISSN 1472-2739 
(on-line) 1472-2747 (printed)\hfill {\lnum\number\count0}\else\ifodd\count0
\def\\{ }\ifx\theshorttitle\relax \thetitle \else\theshorttitle\fi\hfill
{\lnum\number\count0}\else\def\\{ and }{\lnum\number\count0}
\hfill\ifx\theshortauthors\relax 
\theauthors\else\theshortauthors\fi\fi\fi}}\def\@evenhead{\@oddhead}
\def\@oddfoot{\small\lfoot\ifnum\count0=\startpage\copyright\ \gtp\hfill\else
\agt, Volume \thevolumenumber\ (\thevolumeyear)\hfill\fi}
\def\@evenfoot{\@oddfoot}
\makeatother
\fi
\let\maketitlepage\makeagttitle
\let\makeshorttitle\maketitlepage
\let\maketitle\maketitlepage

\newwrite\gtoutfile
\long\gdef\makeheadfile{  
{\def\\{, }\def\s{ }
\immediate\openout\gtoutfile head.xxx
\immediate\write\gtoutfile{Proxy-for: \ifx\theasciiauthors\relax
\theauthors\else\theasciiauthors\fi\s<\ifx\theasciiemail\relax\theemail\else\theasciiemail\fi>}
\immediate\write\gtoutfile{\noexpand\\}
\immediate\write\gtoutfile{Authors: \ifx\theasciiauthors\relax
\theauthors\else\theasciiauthors\fi}
{\def\\{ }\immediate\write\gtoutfile{Title: \ifx\theasciititle\relax
\thetitle\else\theasciititle\fi}}
\immediate\write\gtoutfile{Subj-class: GT or SG, GR etc}
\immediate\write\gtoutfile{MSC-class: \theprimaryclass\ifx\thesecondaryclass\relax\else, \thesecondaryclass\fi}
\immediate\write\gtoutfile{Journal-ref: Algebr. Geom. Topol. \thevolumenumber\s
(\thevolumeyear) \startpage-\finishpage}
\immediate\write\gtoutfile{Comments: Published by Algebraic and
Geometric Topology at}
\immediate\write\gtoutfile{\s\s\s  http://www.maths.warwick.ac.uk/agt/AGTVol\thevolumenumber/agt-\thevolumenumber-\thepapernumber.abs.html}
\immediate\write\gtoutfile{\noexpand\\}
\immediate\write\gtoutfile{}
\ifx\theasciiabstract\relax
\immediate\write\gtoutfile{\theabstract}\else
\immediate\write\gtoutfile{\theasciiabstract}\fi
\immediate\write\gtoutfile{}
\immediate\write\gtoutfile{\noexpand\\}
\immediate\write\gtoutfile{}
\immediate\closeout\gtoutfile}}  

\def\maketitlepage{\makeagttitle\makeheadfile}
\let\makeshorttitle\maketitlepage
\let\maketitle\maketitlepage

\lognumber{44}
\volumenumber{3}
\volumeyear{2003}
\papernumber{44}
\published{20 December 2003}
\pagenumbers{1257}{1273}
\received{20 January 2003}
\revised{12 December 2003}
\accepted{18 December 2003}

\usepackage{amsmath,amssymb}

\newcommand{\ZZ}{\mathbb{Z}}
\newcommand{\CC}{\mathbb{C}}

\newcommand{\PP}{\mathcal{P}}

\newcommand{\disc}{\mathsf{D}}
\newcommand{\lig}{\operatorname{LieGriess}}
\newcommand{\lie}{\operatorname{Lie}}
\newcommand{\griess}{\Sigma\operatorname{Griess}}
\newcommand{\com}{\operatorname{Com}}
\newcommand{\ram}{\operatorname{Ram}}
\newcommand{\alra}{\mathsf{R}}
\newcommand{\odi}{\mathcal{O}_{\mathsf{D}^2}}

\newtheorem{prop}{Proposition}
\newtheorem{lemma}{Lemma}
\newtheorem*{question}{Question}

\begin{document}

\title[On a Hopf Operad]{On a Hopf operad containing the Poisson operad} \date{\today}
\author{Fr\'ed\'eric Chapoton}
\asciiauthors{Frederic Chapoton}
\coverauthors{Fr\noexpand\'ed\noexpand\'eric Chapoton}
\address{LACIM\\
Universit\'e du Qu\'ebec \`a Montr\'eal\\
CP 8888 succ.\ centre ville\\
Montr\'eal Qu\'ebec H3C 3P8\\
Canada}
\asciiaddress{LACIM, Universite du Quebec a
Montreal\\CP 8888 succ. centre ville, Montreal Quebec H3C 3P8, Canada}
\email{chapoton@math.uqam.ca}

\begin{abstract}
  A new Hopf operad $\ram$ is introduced, which contains both the
  well-known Poisson operad and the Bessel operad introduced
  previously by the author. Besides, a structure of cooperad $\alra$
  is introduced on a collection of algebras given by generators and
  relations which have some similarity with the Arnold relations for
  the cohomology of the type $A$ hyperplane arrangement. A map from
  the operad $\ram$ to the dual operad of $\alra$ is defined which we
  conjecture to be a isomorphism.
\end{abstract}

\asciiabstract{%
  A new Hopf operad Ram is introduced, which contains both the
  well-known Poisson operad and the Bessel operad introduced
  previously by the author. Besides, a structure of cooperad R
  is introduced on a collection of algebras given by generators and
  relations which have some similarity with the Arnold relations for
  the cohomology of the type A hyperplane arrangement. A map from
  the operad Ram to the dual operad of R is defined which we
  conjecture to be a isomorphism.}


\maketitle

\setcounter{section}{-1}

\section{Introduction}

The theory of operads has roots in algebraic topology. One well-known
way to build algebraic operads is to start from an operad of
topological spaces and apply the homology functor. A famous example,
due to Cohen \cite{cohen1,cohen2}, is given by the little discs operad
whose homology is the Gerstenhaber operad. The operads defined in this
way inherit more structure from the diagonal of topological spaces:
they are in fact Hopf operads. This phenomenon is similar to the
existence of a bialgebra structure on the homology of a topological
monoid.

This article introduces two algebraic objects. The first one is a Hopf
operad called the Ramanujan operad and denoted by $\ram$, which
contains both the well-known Poisson operad and the Bessel operad
introduced in \cite{hopfbessel}. The second one is a Hopf cooperad
$\alra$, which means that the space $\alra(I)$ associated to a finite
set $I$ is an associative algebra and the cocomposition maps are
morphisms of algebras.

The operad $\ram$ is conjectured to be isomorphic to the linear dual
operad $\alra^*$ of the cooperad $\alra$. A morphism of operad from $\ram$
to $\alra ^*$ is defined, which should give the desired isomorphism.

One motivation for these constructions is an analogy with the case of
the Gerstenhaber operad. The dual algebras of the coalgebras
underlying the Gerstenhaber operad can be presented by generators and
relations, by a theorem of Arnold on the cohomology of the
complexified hyperplane arrangements of type $A$. From this, one can
reach a simple description of the dual cooperad of the Gerstenhaber
operad. This alternative dual description of the Gerstenhaber operad
is sketched at the end of the paper.

There seems to be some kind of similar relation between the cooperad
$\alra$ and some differential forms on the complexified hyperplane
arrangements of type $A$. This relation was already proposed for the
Bessel suboperad in \cite{hopfbessel}. The algebras underlying $\alra$
are given by generators and relations which have some resemblance with
the Arnold relations and seem to contain the relations satisfied by
some simple differential forms.

\smallskip

After some preliminary material on operads in the first section, the
Ramanujan operad is defined in the second section by a distributive
law between the commutative operad and an operad mixing the Lie operad
and the suspended Griess operad. This name has been chosen because the
dimensions are supposed to be given by the so-called Ramanujan
polynomials \cite{chenguo}.

In the next section, the cooperad $\alra$ is defined on a collection
of algebras given by generators and relations. The cocomposition is
motivated by the analogy with the case of the Gerstenhaber operad.
Then a morphism from $\ram$ to $\alra^*$ is defined. Some algebras of
differential forms are introduced, which should be related to $\alra$.
Last, a construction is sketched for the Gerstenhaber operad, which
motivated the formula for cocomposition in $\alra$.

\section{Operads as functors}

Because our language for operads differs in aspect from the most
frequently used setups, this section gathers conveniently some
conventions and definitions.

An operad $\PP$ is a functor from the category of finite sets and
bijections to some monoidal category (sets or vector spaces for
example) together with some extra structure given by composition maps.
Finite sets will be denoted by capital letters $I,J,K$ and so on.
Elements of finite sets will be denoted by letters $i,j,k$ and so on.
In some sense, $i$,$j$,$k$ can be considered as abstract variables when
they are used to denote elements of an arbitrary finite set. The
symbols $\star$ and $\#$ are used as place-holders for composition
maps.

The composition map $\circ_\star$ is defined for any two finite sets $I$
and $J$ as a map from $\PP(I \sqcup \{\star\})\otimes \PP(J) $ to $\PP(I
\sqcup J)$. These composition maps have to satisfy some natural
axioms. Other symbols such as $\#$ are used instead of $\star$ when
iterated compositions appear.

A presentation by generators and relations of an operad is given as
follows: some generators labelled by their inputs, with some specific
symmetry properties with respect to the symmetric group on these
inputs, and some relations involving compositions of these generators.
Consider for example the $\lie$ operad. The generators are $L_{i,j}$
on any set $\{i,j\}$, which stand for the Lie bracket. The generator
$L_{i,j}$ is antisymmetric under the exchange of $i$ and $j$. The
relations are the Jacobi identities (see (\ref{jacobi}) below) on any
set $\{i,j,k\}$, involving generators on various subsets of cardinality
two of $\{i,j,k,\star\}$.

\section{The Ramanujan operad}

In this section, a Hopf operad $\ram$ is defined by a distributive
law. This is similar to the usual definition of the Gerstenhaber
operad by a distributive law between the commutative operad and the
suspended Lie operad.

\subsection{The $\lig$ operad}

The ground field is $\CC$. The ambient category is the monoidal
category of $(\ZZ,\ZZ)$-bigraded vector spaces endowed with two
differentials of respective degree $(1,0)$ and $(-1,0)$. The Koszul
sign rules for the symmetry isomorphisms of the tensor product apply
only with respect to the first degree. The second degree does not play
any role with respect to signs in the formulas.

One can remark that the second degree coincide in the objects
considered here with the eigenvalue of the Laplacian associated to the
pair of opposite differentials.

In this section, an operad $\lig$ is defined which contains the
operads $\lie$ and $\griess$ defining Lie algebras and suspended
commutative non-associative algebras.

The operad $\lig$ is generated by the $\lie$ generator $L_{i,j}$
antisymmetric of degree $(0,1)$ and the $\griess$ generator
$\Omega_{i,j}$ antisymmetric of degree $(1,1)$ modulo the following
relations. 

First, the Jacobi identity defining the Lie operad:
\begin{equation}
  \label{jacobi}
    \sum_{cycl} L_{i,\star}\circ_\star L_{j,k}=0,
\end{equation}
where the summation is over cyclic permutations of $i,j,k$.

Second, a mixed relation between the Lie generator and the $\griess$
generator:
\begin{equation}
  \label{mixte}
  \sum_{cycl} \left( \Omega_{i,\star}\circ_\star L_{j,k} 
  +L_{i,\star}\circ_\star \Omega_{j,k}\right) =0.
\end{equation}
Note that the $\griess$ operad is free on its generator, so there is no
relation involving only $\Omega$.

\subsection{Distributive law}

For the notion of distributive law between operads, see \cite{markl}.

First recall the $\com$ operad, which is generated by $E_{i,j}$
symmetric of degree $(0,0)$ modulo the relation of associativity:
\begin{equation}
  \label{associatif}
  E_{i,\star} \circ_\star E_{j,k}=E_{j,\star} \circ_\star E_{k,i}.
\end{equation}
Then consider the following relations:
\begin{align}
    \label{leibniz}
    L_{i,\star} \circ_{\star} E_{j,k}&=E_{j,\star} \circ_\star
    L_{i,k} +E_{k,\star} \circ_\star
    L_{i,j},\\
    \label{bessel}
    \Omega_{i,\star} \circ_{\star} E_{j,k}&=E_{j,\star} \circ_\star
    \Omega_{i,k} +E_{k,\star} \circ_\star
    \Omega_{i,j}.
\end{align}

\begin{prop}
  The relations (\ref{leibniz}) and (\ref{bessel}) define a
  distributive law from\break
   $\lig \circ \com$ to $\com \circ \lig$. The
  resulting operad is called $\ram$.
\end{prop}

\begin{proof}
  The relation (\ref{leibniz}), which is the Leibniz relation, is
  already known to define a distributive law from $\lie \circ \com$ to
  $\com \circ \lie$. The resulting operad is the Poisson operad.
  
  The relation (\ref{bessel}) is also known to be a distributive law
  from $\griess \circ \com$ to $\com \circ \griess$ which defines the
  Bessel operad, see \cite{hopfbessel}.
  
  So there remains only one condition to check, which comes from
  relation (\ref{mixte}). One has to check that
  \begin{multline}
    (\Omega_{i,\star}\circ_\star L_{j,\#}+\Omega_{j,\star}\circ_\star
    L_{\#,i}+\Omega_{\#,\star}\circ_\star L_{i,j}+\\
    L_{i,\star}\circ_\star \Omega_{j,\#}+L_{j,\star}\circ_\star
    \Omega_{\#,i}+L_{\#,\star}\circ_\star \Omega_{i,j})\circ_{\#} E_{k,\ell},
  \end{multline}
  once rewritten using the distributive laws, reduces to zero modulo
  the relations.

  The result of rewriting is
  \begin{multline*}
    \Omega_{i,\star} \circ_\star 
    (E_{\ell,\#}\circ_\# L_{j,k} +E_{k,\#}\circ_\# L_{j,\ell} )
  - \Omega_{j,\star} \circ_\star 
    (E_{\ell,\#}\circ_\# L_{i,k} +E_{k,\#}\circ_\# L_{i,\ell} )
  \\+ L_{i,\star} \circ_\star 
    (E_{\ell,\#}\circ_\# \Omega_{j,k} +E_{k,\#}\circ_\# \Omega_{j,\ell} )
  -  L_{j,\star} \circ_\star 
    (E_{\ell,\#}\circ_\# \Omega_{i,k} +E_{k,\#}\circ_\# \Omega_{i,\ell} )\\
    -(E_{\ell,\#}\circ_\# \Omega_{\star,k}
    +E_{k,\#}\circ_\# \Omega_{\star,\ell})\circ_\star L_{i,j}
    -(E_{\ell,\#}\circ_\# L_{\star,k}
    +E_{k,\#}\circ_\# L_{\star,\ell})\circ_\star \Omega_{i,j}.
  \end{multline*}

  This becomes, after a second application of the distributive laws,
  \begin{multline*}
    (E_{\#,\star}\circ_\star \Omega_{i,\ell}
    +E_{\ell,\star}\circ_\star \Omega_{i,\#})\circ_\# L_{j,k} 
    +(E_{\#,\star}\circ_\star \Omega_{i,k}
    +E_{k,\star}\circ_\star \Omega_{i,\#})\circ_\# L_{j,\ell} \\
-    (E_{\#,\star}\circ_\star \Omega_{j,\ell}
    +E_{\ell,\star}\circ_\star \Omega_{j,\#})\circ_\# L_{i,k} 
    -(E_{\#,\star}\circ_\star \Omega_{j,k}
    +E_{k,\star}\circ_\star \Omega_{j,\#})\circ_\# L_{i,\ell} \\
+    (E_{\#,\star}\circ_\star L_{i,\ell}
    +E_{\ell,\star}\circ_\star L_{i,\#})\circ_\# \Omega_{j,k} 
    +(E_{\#,\star}\circ_\star L_{i,k}
    +E_{k,\star}\circ_\star L_{i,\#})\circ_\# \Omega_{j,\ell} \\
-    (E_{\#,\star}\circ_\star L_{j,\ell}
    +E_{\ell,\star}\circ_\star L_{j,\#})\circ_\# \Omega_{i,k} 
    -(E_{\#,\star}\circ_\star L_{j,k}
    +E_{k,\star}\circ_\star L_{j,\#})\circ_\# \Omega_{i,\ell} \\
   - E_{\ell,\#} \circ_\# 
   (\Omega_{\star,k} \circ_\star L_{i,j}+L_{\star,k} \circ_\star \Omega_{i,j})
   - E_{k,\#} \circ_\# 
   (\Omega_{\star,\ell} \circ_\star L_{i,j}+L_{\star,\ell}
   \circ_\star \Omega_{i,j}).
  \end{multline*}
  
  Now all $8$ terms starting with $E_{\#,\star}$ annihilates pairwise
  and one can separate in what remains terms starting with
  $E_{\ell,\#}$ and with $E_{k,\#}$. Each of these separate sums is
  zero modulo relation (\ref{mixte}).
\end{proof}

\smallskip

The bigraded dimensions of the $\ram$ operad are yet to be computed.
As explained in the introduction, they should be given by the
Ramanujan polynomials \cite{chenguo}, which are polynomials in
$\{x,y\}$ defined by
\begin{align}
  \psi_1&=1,\\
  \psi_{n+1}&=\psi_n + (x+y) (n \psi_n +x \partial_x \psi_n) \quad n
  \geq 1.
\end{align}
More precisely, the dimension of the homogeneous component of degree
$(i,j)$ of $\ram(\{1,2,\dots,n\})$ should be the coefficient of $x ^i
y ^{j-i}$ in $\psi_n$.

This has been checked for sets with at most four elements. Besides,
the parts of the bigraded dimensions corresponding to the Poisson and
Bessel suboperads are correct, \textit{i.e.}\ match the well-known dimensions of
Poisson and the dimensions of Bessel computed in \cite{hopfbessel}.

\subsection{Hopf structure}

In this section, a coproduct is defined which is compatible with
composition, \textit{i.e.}\ composition becomes a morphism of coalgebras.

The coproduct $\Delta$ is defined on generators by
\begin{equation}
  \label{coprod_delta}
  \begin{cases}
    \Delta(E_{i,j})=E_{i,j} \otimes E_{i,j},\\
    \Delta(L_{i,j})=E_{i,j} \otimes L_{i,j}+L_{i,j} \otimes E_{i,j},\\
    \Delta(\Omega_{i,j})=E_{i,j} \otimes \Omega_{i,j}+
    \Omega_{i,j} \otimes E_{i,j}.
  \end{cases}
\end{equation}

\begin{prop}
  The coproduct $\Delta$ endows $\ram$ with a structure of Hopf operad.
\end{prop}

\begin{proof}
  One has to check the compatibility of $\Delta$ with all relations.
  
  The cases of relations (\ref{jacobi}), (\ref{associatif}) and
  (\ref{leibniz}) are well-know from the Hopf structure of the Poisson
  operad.

  The case of relation (\ref{bessel}) is a consequence of the study of the
  Bessel operad in \cite{hopfbessel}.
  
  So there remains only to check the compatibility of the relation
  (\ref{mixte}). Its coproduct is
  \begin{multline*}
    \sum_{cycl}
    ( E_{i,\star} \circ_{\star} E_{j,k} ) 
    \otimes (\Omega_{i,\star} \circ_{\star} L_{j,k})+
    ( E_{i,\star} \circ_{\star} L_{j,k} ) 
    \otimes (\Omega_{i,\star} \circ_{\star} E_{j,k})\\+
    ( \Omega_{i,\star} \circ_{\star} E_{j,k} ) 
    \otimes (E_{i,\star} \circ_{\star} L_{j,k})+
    ( \Omega_{i,\star} \circ_{\star} L_{j,k} ) 
    \otimes (E_{i,\star} \circ_{\star} E_{j,k})\\+
    ( L_{i,\star} \circ_{\star} \Omega_{j,k} ) 
    \otimes (E_{i,\star} \circ_{\star} E_{j,k})+
    ( L_{i,\star} \circ_{\star} E_{j,k} ) 
    \otimes (E_{i,\star} \circ_{\star} \Omega_{j,k})\\+
    ( E_{i,\star} \circ_{\star} \Omega_{j,k} ) 
    \otimes (L_{i,\star} \circ_{\star} E_{j,k})+
    ( E_{i,\star} \circ_{\star} E_{j,k})
    \otimes (L_{i,\star} \circ_{\star} \Omega_{j,k} ).
  \end{multline*}

  By relations (\ref{associatif}) and (\ref{mixte}), this becomes
  \begin{multline*}
    \sum_{cycl}
    ( E_{i,\star} \circ_{\star} L_{j,k} ) 
    \otimes (\Omega_{i,\star} \circ_{\star} E_{j,k})+
    ( \Omega_{i,\star} \circ_{\star} E_{j,k} ) 
    \otimes (E_{i,\star} \circ_{\star} L_{j,k})\\+
    ( L_{i,\star} \circ_{\star} E_{j,k} ) 
    \otimes (E_{i,\star} \circ_{\star} \Omega_{j,k})+
    ( E_{i,\star} \circ_{\star} \Omega_{j,k} ) 
    \otimes (L_{i,\star} \circ_{\star} E_{j,k}).
  \end{multline*}

  By the distributive laws (\ref{leibniz}) and (\ref{bessel}), this equals
  \begin{multline*}
    \sum_{cycl}
    ( E_{i,\star} \circ_{\star} L_{j,k} ) 
    \otimes (E_{j,\star} \circ_\star
    \Omega_{i,k})+
    ( E_{i,\star} \circ_{\star} L_{j,k} ) 
    \otimes (E_{k,\star} \circ_\star
    \Omega_{i,j})\\+
    (E_{j,\star} \circ_\star
    \Omega_{i,k}) 
    \otimes (E_{i,\star} \circ_{\star} L_{j,k})+
    (E_{k,\star} \circ_\star
    \Omega_{i,j}) 
    \otimes (E_{i,\star} \circ_{\star} L_{j,k})\\+
    (E_{j,\star} \circ_\star
    L_{i,k} ) 
    \otimes (E_{i,\star} \circ_{\star} \Omega_{j,k})+
    (E_{k,\star} \circ_\star
    L_{i,j}) 
    \otimes (E_{i,\star} \circ_{\star} \Omega_{j,k})\\+
    ( E_{i,\star} \circ_{\star} \Omega_{j,k} ) 
    \otimes (E_{j,\star} \circ_\star
    L_{i,k})+
    ( E_{i,\star} \circ_{\star} \Omega_{j,k} ) 
    \otimes (E_{k,\star} \circ_\star
    L_{i,j}).
  \end{multline*}
  
  Consider separately the terms of the form $(E \circ L)\otimes (E
  \circ \Omega)$:
   \begin{multline*}
    \sum_{cycl}
    ( E_{i,\star} \circ_{\star} L_{j,k} ) 
    \otimes (E_{j,\star} \circ_\star
    \Omega_{i,k})+
    ( E_{i,\star} \circ_{\star} L_{j,k} ) 
    \otimes (E_{k,\star} \circ_\star
    \Omega_{i,j})\\+
    (E_{j,\star} \circ_\star
    L_{i,k} ) 
    \otimes (E_{i,\star} \circ_{\star} \Omega_{j,k})+
    (E_{k,\star} \circ_\star
    L_{i,j}) 
    \otimes (E_{i,\star} \circ_{\star} \Omega_{j,k}).
  \end{multline*}

  Shifting the cyclic sum gives
   \begin{multline*}
    \sum_{cycl}
    ( E_{i,\star} \circ_{\star} L_{j,k} ) 
    \otimes (E_{j,\star} \circ_\star
    \Omega_{i,k})+
    ( E_{i,\star} \circ_{\star} L_{j,k} ) 
    \otimes (E_{k,\star} \circ_\star
    \Omega_{i,j})\\+
    (E_{i,\star} \circ_\star
    L_{k,j} ) 
    \otimes (E_{k,\star} \circ_{\star} \Omega_{i,j})+
    (E_{i,\star} \circ_\star
    L_{j,k}) 
    \otimes (E_{j,\star} \circ_{\star} \Omega_{k,i}),
  \end{multline*}
  which is zero. The same is true for terms of the form $(E \circ
  \Omega) \otimes (E \circ L)$. Therefore the coproduct of relation
  (\ref{mixte}) is zero and the proposition is proved.
\end{proof}

\subsection{Two Differentials}

Here are defined two differentials which are derivations for the
composition and coderivations for the coproduct.

The differentials $D$ and $D'$ are defined on generators by
\begin{equation}
  \label{grand_diff}
  \begin{cases}
    D'(E_{i,j})=0,\\
    D'(L_{i,j})=\Omega_{i,j},\\
    D'(\Omega_{i,j})=0.
  \end{cases}
  \begin{cases}
    D(E_{i,j})=0,\\
    D(\Omega_{i,j})=L_{i,j},\\
    D(L_{i,j})=0.
  \end{cases}
\end{equation}

\begin{prop}
  The differentials $D$ and $D'$ can be extended to derivations for
  the composition.
\end{prop}

\begin{proof}
  It is an easy exercise to check against all relations that the
  differentials can be extended to derivations.
\end{proof}

\begin{prop}
  The differentials $D$ and $D'$ are coderivations for the coproduct
  $\Delta$.
\end{prop}

\begin{proof}
  This follows immediately by checking on generators using relations
  (\ref{coprod_delta}) and (\ref{grand_diff}).
\end{proof}

To summarize the results of this section, the $\ram$ operad is a
bigraded Hopf operad endowed with two differentials, which are
derivations and coderivations, \textit{i.e.}\ $\ram$ is a Hopf operad in
the chosen ambient category.

\section{The cooperad $\alra$}

\subsection{Abstract quotient algebras}

Let $I$ be a finite set. Consider the unital commutative associative
algebra $\alra(I)$ generated by elements $a_{i,j}$ antisymmetric of
degree $(0,1)$ and $b_{i,j}$ antisymmetric of degree $(1,1)$ for all
pairs of distinct elements $i,j$ of $I$ modulo the relations
\begin{align}
  \label{rela2}
  a_{i,j}^2&=0,\\
  \label{relaa}
  a_{i,j}a_{j,k}+a_{j,k}a_{k,i}+a_{k,i}a_{i,j}&=0,\\
  \label{relab}
  b_{i,j}a_{j,k}+b_{j,k}a_{k,i}+b_{k,i}a_{i,j}
  +a_{i,j}b_{j,k}+a_{j,k}b_{k,i}+a_{k,i}b_{i,j}&=0,
\end{align}
the relations
\begin{align}
  \label{relabbn}
  a_{i_0,i_1}b_{i_1,i_2}b_{i_2,i_3} \dots b_{i_n,i_0}&=0,\\
  \label{relbbbn}
  b_{i_0,i_1}b_{i_1,i_2}b_{i_2,i_3} \dots b_{i_n,i_0}&=0,
\end{align}
for $n\geq 1$ where $i_0,i_1,i_2,\dots,i_n$ are pairwise different
elements of $I$, and the $12$-terms relations
\begin{multline}
  \label{relbab}
  bab_{i,j,k,\ell}+ bab_{i,k,j,\ell}+bab_{i,j,\ell,k}+ bab_{i,\ell,j,k}+
  bab_{i,k,\ell,j}+ bab_{i,\ell,k,j}\\+bab_{j,i,k,\ell}+ bab_{j,k,i,\ell}+
  bab_{j,i,\ell,k}+ bab_{j,\ell,i,k}+bab_{k,i,j,\ell}+ bab_{k,j,i,\ell}=0,
\end{multline}
\begin{multline}
  \label{relbbb}
  bbb_{i,j,k,\ell}+ bbb_{i,k,j,\ell}+bbb_{i,j,\ell,k}+ bbb_{i,\ell,j,k}+
  bbb_{i,k,\ell,j}+ bbb_{i,\ell,k,j}\\+bbb_{j,i,k,\ell}+ bbb_{j,k,i,\ell}+
  bbb_{j,i,\ell,k}+ bbb_{j,\ell,i,k}+bbb_{k,i,j,\ell}+ bbb_{k,j,i,\ell}=0,
\end{multline}
where $bab_{i,j,k,\ell}=b_{i,j}a_{j,k}b_{k,\ell}$ and
$bbb_{i,j,k,\ell}=b_{i,j}b_{j,k}b_{k,\ell}$ for short.

Remark that the $12$ terms in relations (\ref{relbab}) and
(\ref{relbbb}) correspond to permutations up to reversal.

\smallskip

Note that the subalgebra of elements of first degree $0$
(\textit{i.e.}\  generated by the elements $a_{i,j}$) has already
appeared in the work of Mathieu on the symplectic and Poisson operads
\cite{mathieu} (see also \cite[\S 4.3]{fomikiri} and \cite{gelfvarch}).

\begin{lemma}
  \label{pourfini}
  One has
  \begin{equation}
    \label{nocycle}
    \sq_{i_0,i_1}\sq_{i_1,i_2}\sq_{i_2,i_3} \dots \sq_{i_n,i_0}=0,
  \end{equation}
  for $n\geq 1$ where $i_0,i_1,i_2,\dots,i_n$ are pairwise different
  elements of $I$ and the empty boxes are filled by $a$ and $b$ in an
  arbitrary way.
\end{lemma}
\begin{proof}
  If there is no $a$ at all, equation (\ref{nocycle}) is just equation
  (\ref{relbbbn}). If there is exactly one $a$, then one can use
  commutativity to assume without further restrictions that this $a$
  is the leftmost letter, which gives equation (\ref{relabbn}).
  Therefore one can assume from now on that there are at least two
  $a$.
  
  The proof is by recursion on the length $n$ of the cycle. If $n=1$,
  then the statement is true by relation (\ref{rela2}). Assume that
  $n\geq 2$ and the statement is true for all integers less than $n$.
  
  The proof is now by another recursion on the shortest chain of $b$
  between two $a$.
  
  Assume first that there are two adjacent $a$ in the cycle, say
  $a_{i_0,i_1}a_{i_1,i_2}$. Then one can use relation (\ref{relaa}) to
  replace $a_{i_0,i_1}a_{i_1,i_2}$ by a sum of two terms in the
  product. Each of the two products obtained contains a shorter cycle
  and therefore vanish.
  
  Assume that there are no adjacent $a$ in the cycle. Consider the
  shortest chain of $b$ between two $a$. One can assume without
  restriction that one of these $a$ and one $b$ in the shortest chain
  are $a_{i_0,i_1}b_{i_1,i_2}$. By using relation (\ref{relab}), one
  can replace $a_{i_0,i_1}b_{i_1,i_2}$ by a sum of five terms in the
  product. Among the five products obtained, four have a shorter cycle
  and one has a shorter chain of $b$ between two $a$. Therefore all
  these products vanish by recursion.

  The recursion on the chain is done. The recursion on $n$ is done. 
\end{proof}

\begin{lemma}
  The algebras $\alra(I)$ are finite-dimensional. The second grading
  takes values between $0$ and the cardinality of $I$ minus one.
\end{lemma}
\begin{proof}
  One can map each monomial to a graph on the set $I$ with edges
  colored by $a$ and $b$. If this graph has multiple edges, the
  monomial vanishes by relation (\ref{rela2}) and relations
  (\ref{relabbn}) and (\ref{relbbbn}) for $n=1$. If this graph has a
  loop, the corresponding monomial vanishes by Lemma \ref{pourfini}.
  Therefore only monomials corresponding to forests of simple trees
  can be non-zero in $\alra(I)$. In such a forest, the number of edges
  is at most one less than the cardinality of $I$. As the generators
  have second degree $1$, the maximal second degree of a non-zero
  monomial is therefore bounded by the cardinality of $I$ minus one.
\end{proof}

\begin{lemma}
  \label{diff_aab}
  One has
  \begin{equation}
    \label{aab_24}
    \sum_{\sigma} (aab)_{\sigma}=0,
  \end{equation}
  where $\sigma$ runs over the set of permutations of $\{i,j,k,\ell\}$
  and
  \begin{equation*}
    (aab)_{\sigma}
    =a_{\sigma(i),\sigma(j)}a_{\sigma(j),\sigma(k)}b_{\sigma(k),\sigma(\ell)}.
  \end{equation*}
  One has 
  \begin{equation}
    \sum_{(i',i'')}  T_{i'}^{i''} =0,
  \end{equation}
  $$T_{i}^j=b_{i,j} a_{i,k} a_{i,\ell}\leqno{\hbox{where}}$$ and the sum is over pairs
  of distinct elements of $\{i,j,k,\ell\}$.
\end{lemma}

\begin{proof}
  Let us call $\Sigma$ the first sum and $T$ the second sum.

  Both statements are proved simultaneously. Consider the simplex
  with vertex set $\{i,j,k,\ell\}$. To a facet $f$, one can associate
  a relation $r(f)$ of type (\ref{relab}) and a relation $s(f)$ of
  type (\ref{relaa}). To an edge $e$, one can associate an
  element $a(e)$ and an element $b(e)$.
  
  By summing (with appropriate signs) $r(f)a(e)$ over the set of pairs
  $(f,e)$ where $f$ is a facet and $e$ an edge such that $e \not
  \subset f$, one gets that $\Sigma+2 T$ vanishes.
   
  By summing (with appropriate signs) $s(f)b(e)$ over the set of pairs
  $(f,e)$ where $f$ is a facet and $e$ an edge such that $e \not
  \subset f$, one gets that $\Sigma +T$ vanishes.
\end{proof}

\begin{lemma}
  \label{diff_abb}
  One has
   \begin{equation}
    \label{abb_24}
    \sum_{\sigma} (abb)_{\sigma}=0,
  \end{equation}
  where $\sigma$ runs over the set of permutations of $\{i,j,k,\ell\}$ 
  and
  \begin{equation*}
    (abb)_{\sigma}
    =a_{\sigma(i),\sigma(j)}b_{\sigma(j),\sigma(k)}b_{\sigma(k),\sigma(\ell)}.
  \end{equation*}
\end{lemma}

\begin{proof}
  Consider the simplex with vertex set $\{i,j,k,\ell\}$. To a facet
  $f$, one can associate a relation $r(f)$ of type (\ref{relab}). To
  an edge $e$, one can associate an element $b(e)$.

  By summing (with appropriate signs) $r(f)b(e)$ over the set of pairs
  $(f,e)$ where $f$ is a facet and $e$ an edge such that $e \not
  \subset f$, one gets that the sum (\ref{abb_24}) vanishes.
\end{proof}

One can define two differentials $d$ and $d'$ on generators by 
\begin{equation}
  \label{petit_diff}
  \begin{cases}
    d(a_{i,j})=b_{i,j},\\
    d(b_{i,j})=0.
  \end{cases}
  \begin{cases}
    d'(b_{i,j})=a_{i,j},\\
    d'(a_{i,j})=0.
  \end{cases}
\end{equation}
\eject

\begin{prop}
  The differentials $d$ and $d'$ can be extended to derivations of the
  algebra $\alra(I)$.
\end{prop}

\begin{proof}
  The check is quite easy for the differential $d$. For $d'$, the
  only non-trivial cases are relations (\ref{relabbn}), (\ref{relbab}) and
  (\ref{relbbb}). 

  The case of relation (\ref{relabbn}) is settled by Lemma \ref{pourfini}.

  The image by $d'$ of relation (\ref{relbab}) is
  exactly the sum (\ref{aab_24}) which vanishes by Lemma
  \ref{diff_aab}.
  
  The image by $d'$ of relation (\ref{relbbb}) is the sum of relation
  (\ref{relbab}) and relation (\ref{abb_24}) and therefore vanishes by
  Lemma \ref{diff_abb}.
\end{proof}

\subsection{Cocomposition maps}

Let $I$ and $J$ be two finite sets. Motivated by the similar
cocomposition (\ref{coco_gerst}) for the dual of the Gerstenhaber
operad, one defines the cocomposition map $\Theta^{\star}_{I,J}$ from
$I \sqcup J$ to $(I \sqcup \{\star\}, J)$ on generators by
\begin{equation}
  \Theta^{\star}_{I,J}(a_{i,j})=
  \begin{cases}
    a_{i,j} \otimes 1\text{ if }i,j \in I,\\
    1 \otimes a_{i,j}\text{ if }i,j \in J,\\
    a_{i, \star} \otimes 1\text{ if }i \in I\text{ and }j \in J
  \end{cases}
\end{equation}
and
\begin{equation}
  \Theta^{\star}_{I,J}(b_{i,j})=
  \begin{cases}
    b_{i,j} \otimes 1\text{ if }i,j \in I,\\
    1 \otimes b_{i,j}\text{ if }i,j \in J,\\
    b_{i, \star} \otimes 1\text{ if }i \in I\text{ and }j \in J.
  \end{cases}
\end{equation}

\begin{prop}
  This defines morphisms $\Theta^{\star}_{I,J}$ of bidifferential
  algebras from $\alra(I\sqcup J)$ to $\alra(I\sqcup \{\star\})
  \otimes \alra(J)$.
\end{prop}
\begin{proof}
  First one has to check against all relations that
  $\Theta^{\star}_{I,J}$ can be extended to a morphism of algebras. By
  the very simple shape of cocomposition, the compatibility is clear
  if all indices involved are in $J$ or if all but maybe one are in
  $I$.
  
  So one can assume that there is at least one index in $I$ and at
  least two indices in $J$. Again compatibility is easy to check for
  all relations involving at most three indices. The only non-trivial
  cases are the relations (\ref{relabbn}) and (\ref{relbbbn}) and the
  12-terms relations (\ref{relbab}) and (\ref{relbbb}).
  
  Consider first the case of relations (\ref{relabbn}) and
  (\ref{relbbbn}). More generally, consider any cycle
  $\sq_{i_0,i_1}\sq_{i_1,i_2}\dots\sq_{i_n,i_0}$, where boxes are
  either $a$ or $b$.
  
  As there is at least one index of the cycle in $I$ and at least two
  in $J$, one can assume without restriction that $i_0 \in J$,
  $i_1,\dots,i_k \in I$ and $i_{k+1} \in J$ for some $k\geq 1$. Then
  the left tensor in the image by $\Theta^{\star}_{I,J}$ of the cycle
  contains the cycle
  $\sq_{\star,i_1}\sq_{i_1,i_2}\dots\sq_{i_k,\star}$ for some $a$ and
  $b$ in the boxes, and therefore vanishes by Lemma \ref{pourfini}.

  Now consider for example the case of (\ref{relbab}) with $i,j \in I$
  and $k, \ell \in J$. Its cocomposition is given by
  \begin{multline}
    (b_{i,j}\otimes 1)(a_{j,\star}\otimes 1)(1 \otimes b_{k,\ell})
    +(b_{i,\star}\otimes 1)(a_{\star,j}\otimes 1)(b_{j,\star}\otimes 1)\\
    +(b_{i,j}\otimes 1)(a_{j,\star}\otimes 1)(1 \otimes b_{\ell,k})
    +(b_{i,\star}\otimes 1)(a_{\star,j}\otimes 1)(b_{j,\star}\otimes 1)\\
    +(b_{i,\star}\otimes 1)(1 \otimes a_{k,\ell})(b_{\star,j}\otimes
    1)
    +(b_{i,\star}\otimes 1)(1 \otimes a_{\ell,k})(b_{\star,j}\otimes 1)\\
    +(b_{j,i}\otimes 1)(a_{i,\star}\otimes 1)(1 \otimes b_{k,\ell})
    +(b_{j,\star}\otimes 1)(a_{\star,i}\otimes 1)(b_{i,\star}\otimes 1)\\
    +(b_{j,i}\otimes 1)(a_{i,\star}\otimes 1)(1 \otimes b_{\ell,k})
    +(b_{j,\star}\otimes 1)(a_{\star,i}\otimes 1)(b_{i,\star}\otimes 1)\\
    +(b_{\star,i}\otimes 1)(a_{i,j}\otimes 1)(b_{j,\star}\otimes 1)
    +(b_{\star,j}\otimes 1)(a_{j,i}\otimes 1)(b_{i,\star}\otimes 1).
  \end{multline}
  This is equal to
  \begin{multline}
    (b_{i,j} a_{j,\star}) \otimes b_{k,\ell} +(b_{i,\star} a_{\star,j}
    b_{j,\star}) \otimes 1 +(b_{i,j} a_{j,\star}) \otimes b_{\ell,k}
    +(b_{i,\star}a_{\star,j}b_{j,\star})\otimes 1\\
    +(b_{i,\star}b_{\star,j}) \otimes a_{k,\ell}
    +(b_{i,\star}b_{\star,j}) \otimes a_{\ell,k} +(b_{j,i}a_{i,\star})
    \otimes b_{k,\ell}
    +(b_{j,\star}a_{\star,i}b_{i,\star}) \otimes 1\\
    +(b_{j,i}a_{i,\star}) \otimes b_{\ell,k}
    +(b_{j,\star}a_{\star,i}b_{i,\star})\otimes 1
    +(b_{\star,i}a_{i,j}b_{j,\star})\otimes 1
    +(b_{\star,j}a_{j,i}b_{i,\star})\otimes 1.
  \end{multline}
  Using antisymmetry and some relations, this is seen to be zero. The
  proof in the remaining cases for (\ref{relbab}) and (\ref{relbbb})
  is similar and left to the reader.
  
  This map clearly respects both differentials, as can be checked on
  generators.
\end{proof}

\begin{prop}
  The applications $\Theta$ define a cooperad structure on $\alra$.
\end{prop}
\begin{proof}
  One has to check on the generators of $\alra(I \sqcup J \sqcup K)$
  that
  \begin{equation}
    \label{coassocI}
    (\Theta^{\star}_{I,J\sqcup \{\#\}}\otimes Id_K) \circ 
    \Theta ^{\#}_{I \sqcup J,K}=(Id_{I \sqcup \{\star \}}
    \otimes \Theta ^{\#}_{J,K} ) \circ \Theta ^{\star}_{I, J \sqcup K}
  \end{equation}
  and that 
  \begin{equation}
    (\Theta ^{\star}_{I \sqcup \{\# \},J }\otimes Id_K) \circ 
    \Theta ^{\#}_{I \sqcup J,K}=(Id_I \otimes \tau) \circ 
    (\Theta ^{\#}_{I \sqcup \{\star\},K }\otimes Id_J) \circ 
    \Theta ^{\star}_{I \sqcup K,J},
  \end{equation}
  where $\tau$ is the symmetry isomorphism for $\alra(J)\otimes \alra(K)$. 
  
  The proof is case by case according to the indices of the generator.
  Consider for example equation (\ref{coassocI}) and a generator
  $a_{i,j}$ in $\alra(I \sqcup J \sqcup K)$ with $i \in I$ and $j \in
  K$. One has on the one hand
  \begin{equation*}
    (\Theta^{\star}_{I,J\sqcup \{\#\}}\otimes Id_K) \circ 
    \Theta ^{\#}_{I \sqcup J,K}(a_{i,j})
    =(\Theta^{\star}_{I,J\sqcup \{\#\}}\otimes Id_K)(a_{i,\#} \otimes 1)
      =a_{i,\star} \otimes 1 \otimes 1.
  \end{equation*}
  On the other hand,
  \begin{equation*}
    (Id_{I \sqcup \{\star \}}
    \otimes \Theta ^{\#}_{J,K} ) \circ \Theta ^{\star}_{I, J \sqcup K}
    (a_{i,j})=(Id_{I \sqcup \{\star \}}
    \otimes \Theta ^{\#}_{J,K} )
    (a_{i,\star} \otimes 1)=a_{i,\star} \otimes 1 \otimes 1.
  \end{equation*}
The remaining cases are similar and left to the reader.
\end{proof}

\subsection{Morphism of operads}

Here is defined a morphism $\rho$ from the operad $\ram$ to the dual
operad $\alra^*$ of the cooperad $\alra$.

Consider the dual vector space $\alra^*(I)$ of $\alra(I)$. This vector
space is bigraded. Define elements $1^*$, $b_{i,j}^*$, $a_{i,j}^*$ in
$\alra^*(I)$ as the dual basis (with respect to the pairing $\alra(I)
\otimes \alra ^*(I) \to \CC$) for the components of degree $(0,0)$,
$(0,1)$ and $(1,1)$ respectively.

The map $\rho$ is defined on the generators of $\ram$ by
\begin{equation}
  \begin{cases}
    E_{i,j} \mapsto 1^*,\\
    \Omega_{i,j} \mapsto b_{i,j}^*,\\
    L_{i,j} \mapsto a_{i,j}^*.
  \end{cases}
\end{equation}

\begin{prop}
  This defines a map $\rho$ of Hopf operads from $\ram$ to $\alra^*$. The
  map $\rho$ intertwines $d$ with $D^*$ and $d'$ with $(D')^*$.
\end{prop}

\begin{proof}
  First, one has to check that this indeed defines a morphism of
  operads, \textit{i.e.}\ the compatibility with relations defining
  $\ram$.
  
  For example, let us check the compatibility for relation
  (\ref{mixte}). By the bigrading, it is sufficient to prove that the
  corresponding linear form vanishes on the dual bihomogeneous
  component. First compute the following cocompositions:
  \begin{align}
    \Theta ^{\star}_{\{i\},\{j,k\}}(a_{i,j}b_{j,k})&
    =a_{i,\star} \otimes b_{j,k},\\
    \Theta ^{\star}_{\{i\},\{j,k\}}(a_{j,k}b_{k,i})&
    =b_{\star,i} \otimes a_{j,k},\\
    \Theta ^{\star}_{\{i\},\{j,k\}}(a_{k,i}b_{i,j})&=0,\\
    \Theta ^{\star}_{\{i\},\{j,k\}}(b_{i,j}a_{j,k})&=
    b_{i,\star} \otimes a_{j,k},\\
    \Theta ^{\star}_{\{i\},\{j,k\}}(b_{j,k}a_{k,i})&=
    a_{\star,i} \otimes b_{j,k},\\
    \Theta ^{\star}_{\{i\},\{j,k\}}(b_{k,i}a_{i,j})&=0.
  \end{align}
  From this, one can deduce a description of the linear forms
  $a^*_{i,\star} \circ_{\star} b ^*_{j,k}$ and $b^*_{i,\star}
  \circ_{\star} a ^*_{j,k}$ by their values on a basis of the
  homogeneous component of degree $(1,2)$ of $\alra(\{i,j,k\})$.

  Now the sum (\ref{mixte}) is mapped by $\rho$ to
  \begin{equation}
    \sum_{cycl} \left( b ^*_{i,\star}\circ_\star a ^*_{j,k} 
    +a ^*_{i,\star}\circ_\star b ^*_{j,k} \right).
  \end{equation}
  One then checks that this sum vanishes as a linear form.

  The proof of compatibility for the other relations is similar.
  
  The intertwining property for differentials is clear on the
  generators $\ram(\{i,j\})$ of $\ram$. It is also easy to prove that
  this map is a morphism of coalgebras by checking on generators of
  $\ram$.
\end{proof}

That the map $\rho$ is an isomorphism has been checked for sets with
at most three elements. Furthermore the bigraded dimensions of $\ram$
and $\alra^*$ coincide for sets with at most four elements. One can
therefore ask the following
\begin{question}
  Is $\rho$ an isomorphism ?
\end{question}

Remark that it follows from the fact that $\rho$ is a morphism of Hopf
operads that, by transposition, the relations of the algebras
underlying $\alra$ are satisfied in the dual algebras of the
coalgebras underlying $\ram$.

\subsection{Algebras of differential forms}

Here, a tentative relation of $\alra$ with differential forms on
hyperplane arrangements is proposed.

Let $I$ be a finite set and $\CC^I$ be the vector space with
coordinates $(x_i)_{i \in I}$. Let $\textsf{H}_I$ be the union of all
hyperplanes $x_i-x_j=0$ for $i\not= j$ in the subspace $\sum_{i \in I}
x_i=0$ of $\CC^I$ (this is a type $A$ hyperplane arrangement).
Consider the subalgebra of the algebra of differential forms with
poles along $\textsf{H}_I$ generated over $\CC$ by elements
$a_{i,j}=1/(x_i-x_j)$ and $b_{i,j}=d(1/(x_i-x_j))$ for $i \not=j$
(here $d$ is the de Rham differential).

The elements $a_{i,j}$ and $b_{i,j}$ are antisymmetric. There are two
natural gradings on this algebra: the first one is by the degree as a
differential form, the second one is the homogeneity degree where all
variables $x_i$ are taken homogeneous of degree minus one.

One has the following relations:
\begin{align*}
    a_{i,j}a_{j,k}+a_{j,k}a_{k,i}+a_{k,i}a_{i,j}&=0,\\
    b_{i,j}a_{j,k}+b_{j,k}a_{k,i}+b_{k,i}a_{i,j}
    +a_{i,j}b_{j,k}+a_{j,k}b_{k,i}+a_{k,i}b_{i,j}&=0,
\end{align*}
and
\begin{equation*}
  b_{i_0,i_1}b_{i_1,i_2}b_{i_2,i_3} \dots b_{i_n,i_0}=0,
\end{equation*}
for $n\geq 1$ where $i_0,i_1,i_2,\dots,i_n$ are pairwise different
elements of $I$.

It is plausible that the abstract algebras $\alra(I)$ introduced
before are quotients of these concrete algebras of differential forms.
The main problem is to find some geometric reason for the relations of
$\alra(I)$.

\subsection{On the Gerstenhaber operad}

This section is mainly for motivation and details are therefore
omitted.

Recall the topological little discs operad $\disc^2$, where
$\disc^2(I)$ is the space of disjoint embeddings of scaled unit discs,
bijectively labeled by $I$, inside a unit disc. The composition inside
a little disc is obtained by replacing this little disc by a
collection of little discs appropriately scaled, see \cite{voronov}
for further details.

The algebra $\odi(I)$ defined by
\begin{equation}
  \CC[x_i][(x_i-x_j)^{-1}][[\epsilon_i]][\epsilon_i ^{-1}]
\end{equation}
is an algebraic analog of the algebra of functions on the space
$\disc^2(I)$, where the $x$ variables are the pairwise-different
complex coordinates of the centers and the $\epsilon$ variables are
the infinitesimal non-vanishing real radiuses. Assuming that radiuses
are infinitesimal ensures disjointness of discs. One can
easily translate the composition rule of the topological little discs
operad into a cocomposition rule defining a cooperad on the collection
of the algebras $\odi(I)$ for all finite sets $I$.

Now the Gerstenhaber operad can be defined as the homology of the
little discs operad \cite{gerstvoro,cohen1,cohen2}. As the space
$\disc^2(I)$ is homotopy equivalent to the complement of a
complexified hyperplane arrangement of type $A$, a theorem of Arnold
\cite{arnold} implies that its cohomology is generated by the classes
of the differential forms
\begin{equation}
  \omega_{i,j}=d(\log(x_i-x_j)),
\end{equation}
subject only to the relations
\begin{equation}
  \label{arnold}
  \omega_{ij}\omega_{jk}+\omega_{jk}\omega_{ki}+\omega_{ki}\omega_{ij}=0.
\end{equation}
One can extend the algebraic cocomposition rules for the collection of
algebras $\odi(I)$ obtained before to cocomposition rules for the
collection of algebras of differential forms of $\odi(I)$ with respect
to the $x$ variables. It is then possible to restrict these rules to
the collection of subalgebras generated by the forms $\omega_{i,j}$.
The result is as follows for the cocomposition map
$\Theta^{\star}_{I,J}$ from $I \sqcup J$ to $(I \sqcup \{\star \} ,
J)$:
\begin{equation}
  \label{coco_gerst}
    \Theta^{\star}_{I,J}(\omega_{i,j})=
  \begin{cases}
    \omega_{i,j} \otimes 1\text{ if }i,j \in I,\\
    1 \otimes \omega_{i,j}\text{ if }i,j \in J,\\
    \omega_{i, \star} \otimes 1\text{ if }i \in I\text{ and }j \in J.
  \end{cases}
\end{equation}
Together with the Arnold relations (\ref{arnold}), this provides an
algebraic description of the dual cooperad of the Gerstenhaber operad.

\Addresses\recd


\begin{thebibliography}

\bibitem{arnold}
\textbf{V\,I Arnold}, \emph{The cohomology ring of the group of dyed braids},
  Mat. Zametki 5 (1969) 227--231

\bibitem{hopfbessel}
\textbf{F Chapoton}, \emph{A {H}opf operad of forests of binary trees and
  related finite-dimensional algebras.} (Sept. 2002), {\tt
  arXiv:math.CO/0209038}

\bibitem{chenguo}
\textbf{W\,Y\,C Chen}, \textbf{V\,J\,W Guo}, \emph{Bijections behind the
  {R}amanujan polynomials}, Adv. in Appl. Math. 27 (2001) 336--356, special
  issue in honor of Dominique Foata's 65th birthday (Philadelphia, PA, 2000)

\bibitem{cohen1}
\textbf{F\,R Cohen}, \emph{The homology of ${C}_{n+1}$-spaces, $n\geq 0$},
  from: ``The homology of iterated loop spaces.'', Lecture Notes in Mathematics
  533, Springer (1976)

\bibitem{cohen2}
\textbf{F\,R Cohen}, \emph{Artin's braid groups, classical homotopy theory, and
  sundry other curiosities}, from: ``Braids (Santa Cruz, CA, 1986)'', Contemp.
  Math. 78, Amer. Math. Soc., Providence, RI (1988)  167--206

\bibitem{fomikiri}
\textbf{S Fomin}, \textbf{A\,N Kirillov}, \emph{Quadratic algebras, {D}unkl
  elements, and {S}chubert calculus}, from: ``Advances in geometry'', Progr.
  Math. 172, Birkh\"auser Boston, Boston, MA (1999)  147--182

\bibitem{gerstvoro}
\textbf{M Gerstenhaber}, \textbf{A\,A Voronov}, \emph{Homotopy {$G$}-algebras
  and moduli space operad}, Internat. Math. Res. Notices  (1995) 141--153
  (electronic)

\bibitem{markl}
\textbf{M Markl}, \emph{Distributive laws and {K}oszulness}, Ann. Inst. Fourier
  (Grenoble) 46 (1996) 307--323

\bibitem{mathieu}
\textbf{O Mathieu}, \emph{The symplectic operad}, from: ``Functional analysis
  on the eve of the 21st century, Vol.\ 1 (New Brunswick, NJ, 1993)'', Progr.
  Math. 131, Birkh\"auser Boston, Boston, MA (1995)  223--243

\bibitem{gelfvarch}
\textbf{A Varchenko}, \textbf{I Gelfand}, \emph{Heaviside functions of a
  configuration of hyperplanes}, Funktsional. Anal. i Prilozhen. 21 (1987)
  1--18, 96

\bibitem{voronov}
\textbf{A\,A Voronov}, \emph{The {S}wiss-cheese operad}, from: ``Homotopy
  invariant algebraic structures (Baltimore, MD, 1998)'', Contemp. Math. 239,
  Amer. Math. Soc., Providence, RI (1999)  365--373

\end{thebibliography}
\end{document}